\documentclass[a4paper]{elsarticle}

\usepackage{lineno,hyperref}
\modulolinenumbers[5]

\journal{Journal of \LaTeX\ Templates}

\usepackage{eurosym}
\usepackage{amsmath}
\usepackage{amsthm}
\usepackage{amssymb}
\usepackage{epsfig}
\usepackage{psfrag}
\usepackage{graphicx}
\usepackage[mathcal]{euscript}
\usepackage{color}
\usepackage{helvet}
\usepackage{enumitem}
\usepackage{setspace}

\usepackage{todonotes}

\theoremstyle{plain}
\newtheorem{thm}{Theorem}[section]

\newtheorem{lem}[thm]{Lemma}

\theoremstyle{definition}

\theoremstyle{remark}

\def\P{{\cal P}}

\begin{document}

\begin{frontmatter}

\title{High-order Virtual Element Method \\ on polyhedral meshes}

\author[add1,add2]{L. Beir\~{a}o da Veiga}
\ead{lourenco.beirao@unimib.it}
\author[add1]{F. Dassi}
\ead{franco.dassi@unimib.it}
\author[add1,add2]{A. Russo}
\ead{alessandro.russo@unimib.it}

\address[add1]{Department of Mathematics and Applications, University of Milano - Bicocca,\\ 
Via Cozzi 53, I-20153, Milano (Italy)}
\address[add2]{IMATI-CNR, 27100 Pavia (Italy)}

\begin{abstract}
We develop a numerical assessment of the Virtual Element Method for the discretization of a diffusion-reaction model problem, for higher ``polynomial'' order $k$ and three space dimensions. Although the main focus of the present study is to illustrate some $h$-convergence tests for different orders $k$, we also hint on other interesting aspects such as structured polyhedral Voronoi meshing, robustness in the presence of irregular grids, sensibility to the stabilization parameter and convergence with respect to the order $k$.
\end{abstract}


\end{frontmatter}




\section{Introduction}



The Virtual Element Method (VEM) was introduced in \cite{volley,Hitchhikers} as a generalization of the Finite Element Method (FEM) that allows for very general polygonal and polyhedral meshes, also including non convex and very distorted elements. The VEM is not based on the explicit construction and evaluation of the basis functions, as standard FEM, but on a wise choice and use of the degrees of freedom in order to compute the operators involved in the discretization of the problem. The adopted basis functions are virtual, in the sense that they follow a rigorous definition, include (but are not restricted to) standard polynomials but are \emph{not} computed in practice; the accuracy of the method is guaranteed by the polynomial part of the virtual space. Using such approach introduces other potential advantages, such as exact satisfaction of linear constraints \cite{Stokes:divfree} and the possibility to build easily discrete spaces of high global regularity \cite{Brezzi-Marini:2012,CH-VEM}.
Since its introduction, the VEM has shared a good degree of success and was applied to a large array of problems. We here mention, in addition to the ones above, a sample of papers \cite{Ahmed-et-al:2013,variable-primale,largedefs,nonconforming,Benedetto-VEM-2,Gatica-1,Steklov-VEM,Topology-VEM,vacca2016virtual} and refer to 
\cite{Arxiv-hp-corner}
for a more complete survey of the existing VEM literature.  

Although the construction of the Virtual Element Method for three dimensional problems is accomplished in many papers, 
at the current level of development very few 3D numerical experiments are available in the literature \cite{Paulino-VEM,Topology-VEM,largedefs}. Moreover, all these tests are limited to the lowest order case ($k=1$).

The objective of this work is to numerically validate,  for the first time, the VEM of general order $k$ for three dimensional problems and show that this technology is practically viable also in this case. We consider a simple diffusion-reaction model problem in primal form and follow faithfully the construction in \cite{volley,Ahmed-et-al:2013} and the coding guidelines of  \cite{Hitchhikers}. 
Although the main focus of the present study is to illustrate some standard $h$-convergence tests for different orders $k$, it also hints on other aspects. In particular, we show some interesting possibilities related to polyhedral Voronoi meshing, we underline the robustness of the method in the presence of irregular grids, we investigate its sensibility to the stabilization parameter and consider also a convergence analysis in terms of the ``polynomial'' order $k$.

The paper is organized as follows. In Section \ref{theVEM} we introduce the model problem and the Virtual Element Method in three dimensions. The review of the method is complete but brief, and we refer to other contributions in the literature for a more detailed presentation of the scheme. Afterwards, in Section \ref{NUMS} an array of numerical tests are shown. 


\section{The Virtual Element discretization}\label{theVEM}

In the present section we give a brief overview of the Virtual Element method in three space dimensions for the simple model problem of diffusion-reaction in primal form. More details on the method for this same model and formulation can be found in \cite{volley,Hitchhikers,Ahmed-et-al:2013} while extension to variable coefficients 
is presented in \cite{variable-primale}. In the following $k$ will denote a positive integer number, associated to the ``polynomial degree'' of the virtual element scheme.
 
\subsection{Notation} 

\newcommand{\myboldmath}[1]{{\boldsymbol{#1}}}

\newcommand{\x}{\myboldmath{x}}
\newcommand{\ba}{\myboldmath{\alpha}}
\newcommand{\ma}{m_\ba}

\newcommand{\D}{{D}}

\newcommand{\Dhat}{\widehat{\D}}
\newcommand{\xD}{\x_\D}
\newcommand{\hP}{h_P}
\newcommand{\R}{{\mathbb{R}}}

In the following, $E$ will denote a polygon and
$P$ a polyhedron, while faces, edges and vertices will be 
indicated by $f$, $e$, and $\nu$ respectively.

If $P$ is a polyhedron in $\R^3$, we will denote by $\boldsymbol{x}_P$, $h_P$ and $|P|$
the centroid, the diameter, and the volume of $P$, respectively.
The set of polynomials of degree less than or equal to $s$ in $P$
will be indicated by $\P_s(P)$.
If $\ba=(\alpha_1,\alpha_2,\alpha_3)$ is a multiindex, we will indicate by
$\ma$ the scaled monomial
\begin{equation}
\ma = 
\left(\dfrac{x-x_P}{\hP}\right)^{\alpha_1}
\left(\dfrac{y-y_P}{\hP}\right)^{\alpha_2}
\left(\dfrac{z-z_P}{\hP}\right)^{\alpha_3}
\end{equation}
and we will denote by $\P^{{\rm hom}}_{s}(P)$ the
the space of scaled monomials of degree exactly (and no less than) $s$:
\begin{equation}
\P^{{\rm hom}}_{s}(P) = \text{span} \{\ma \text{ , } |\ba|=s\}
\end{equation}
where $|\ba|=\alpha_1+\alpha_2+\alpha_3$. 
The case of a polygon $E\subset\R^2$ is completely analogous.

A face $f$ of a polyhedron is treated as a two-dimensional set, using local
coordinates $(x,y)$ on the face. Edges of polyhedra and polygons are treated
in an analogous way as one-dimensional set.

\subsection{The model problem}
  
Let $\Omega \subset {\mathbb R}^3$ represent the domain of interest (that we assume to be a polyhedron) and let $\Gamma$ denote a subset of its boundary, that we assume for simplicity to be given by a union of some of its faces. We denote by $\Gamma' = \partial\Omega / \Gamma$. 

We consider the simple reaction-diffusion problem
\begin{equation}
\left\{
\begin{aligned}
-\Delta u + u &= f\quad\text{in }\Omega\\
 u &= r\quad\text{on }\Gamma \\
 \dfrac{\partial u}{\partial n} &= g \quad\text{on }\Gamma'
\end{aligned}
\right.
\end{equation}
where $f  \in L^2(\Omega)$ and $g \in L^2(\Gamma')$ denote respectively the applied load and Neumann boundary data,  
and $r \in H^{1/2}(\Gamma)$ is the assigned boundary data function.
The variational form of our model problem reads
\begin{equation}
\left\{
\begin{aligned}
& \textrm{Find } u \in H^1_\Gamma(\Omega) \textrm{ such that } \\
& \int_\Omega \nabla u \cdot \nabla v  \: + \: \int_\Omega u v \: = \: \int_\Omega f v + \int_{\Gamma'} g v
\qquad \forall v \in H^1_{\Gamma,0}(\Omega) \ ,
\end{aligned}
\right.
\end{equation}
where 
$$
H^1_\Gamma = \big\{ v \in H^1(\Omega) : v|_\Gamma = r \big\} \ , \quad
H^1_{\Gamma,0} = \big\{ v \in H^1(\Omega) : v|_\Gamma = 0 \big\}.
$$

\subsection{Virtual elements on polygons}

We start by defining the virtual element space on polygons. 
Given a generic polygon $E$, let the preliminary virtual space
$$
\widetilde{V}^k(E) = \Big\{ v \in H^1(E) \cap C^0(E) \: : \: v|_e \in \P_k(e)\ \forall e \in \partial E, \ 
\Delta v \in \P_{k}(E) \Big\} 
$$
with $e$ denoting a generic edge of the polygon.

For any edge $e$, let the points $\{ \nu_e^i \}_{i=1}^{k-1}$ be given by the $k-1$ internal points of the Gauss-Lobatto integration rule of order $k+1$ on the edge.
We now introduce three sets of linear operators from $\widetilde{V}^k(E)$ into real numbers.
For all $v$ in $\widetilde{V}^k(E)$:
\begin{align}
\label{dof1}
&\bullet\quad\textrm{evaluation of } v(\nu) \ \forall \nu \textrm{ vertex of } E ; \\
\label{dof2}
&\bullet\quad\textrm{evaluation of } v(\nu_e^i) \ \forall e \in \partial E, \ i=\{1,2,..,k-1 \} ; \\
\label{dof3}
&\bullet\quad\textrm{moments } \int_E v \: p_{k-2} \ \ \forall p_{k-2} \in \P_{k-2}(E) .
\end{align}
The following projector operator $\Pi^\nabla_E \: : \: \widetilde{V}^k(E) \rightarrow \P_k(E)$ will be useful in the definition of our space and also for computational purposes.
For any $v \in \widetilde{V}^k(E)$,  the polynomial $\Pi^\nabla_E  v \in \P_k(E)$ is defined by (see \cite{volley})
$$
\left\{
\begin{aligned}
& \int_E \nabla (v - \Pi^\nabla_E  v)  \cdot \nabla p_k = 0 \quad \forall p_k \in \P_k(E) \\
& \text{for $k=1$:}\quad\sum_{\substack{\nu\,\text{vertex} \\[0.5mm]\text{ of}\,E}} \left(v(\nu) - \Pi^\nabla_E  v(\nu)\right) = 0, \\
& \text{for $k\ge2$:}\quad\int_{ E} (v - \Pi^\nabla_E  v) = 0 .
\end{aligned}
\right.
$$ 
Note that, given any $v \in \widetilde{V}^k(E)$ the polynomial $\Pi^\nabla_E  v$ only depends on the values of the operators \eqref{dof1}-\eqref{dof2}-\eqref{dof3}. Indeed, an integration by parts easily shows that the values of the above operators applied to $v$ are sufficient to uniquely determine $\Pi^\nabla_E  v$ and no other information on the function $v$ is required \cite{volley,Hitchhikers}.

We are now ready to present the two dimensional virtual space:
$$
V^k(E) = \Big\{ v \in \widetilde{V}^k(E) \: : \: \int_E v \, q = \int_E (\Pi^\nabla_E  v) q \textrm{ for all } q \in 
\P_{k-1}^{{\rm hom}}(E)\cup\P_{k}^{{\rm hom}}(E) \Big\} .
$$
It is immediate to check that $\P_k(E) \subseteq V^k(E) \subseteq \widetilde{V}^k(E)$.
Moreover the following lemma holds (the proof can be found in \cite{Ahmed-et-al:2013,variable-primale}). 
\begin{lem}\label{lem:dofs1}
The operators \eqref{dof1}-\eqref{dof2}-\eqref{dof3} constitute a set of degrees of freedom for the space $V^k(E)$.
\end{lem}
Finally note that, since any $p_{k} \in \P_{k}(E)$ can be written in an unique way as $p_{k} = p_{k-2} + q$, with $p_{k-2} \in \P_{k-2}(E)$ and 
$q\in\P_{k-1}^{{\rm hom}}(E)\cup\P_{k}^{{\rm hom}}(E)$, it holds
\begin{equation}\label{obs-1}
\int_E v \, p_{k} = \int_E v \, (p_{k-2} + q)  = 
\int_E v \, p_{k-2} + \int_E (\Pi^\nabla_E  v) q 
\end{equation}
for any $p_{k} \in \P_{k}(E)$. The first term in the right hand side above can be calculated recalling \eqref{dof3} while the second one can be computed directly by integration. This shows that we can actually compute $\int_E v \, p_{k}$ for any 
$p_{k} \in \P_{k}(E)$ by using only information on the degree of freedom values of $v$.

\subsection{Virtual elements on polyhedrons}

Let $\Omega_h$ be a partition of $\Omega$ into non-overlapping and conforming polyhedrons. We start by defining the virtual space $V^k$ locally, on each polyhedron $P \in \Omega_h$. Note that each face $f \in \partial P$ is a two-dimensional polygon.
Let the following boundary space
\newcommand{\BdPk}{{\cal B}^k({\partial P})}
$$
\BdPk = \Big\{ v \in C^0(\partial P) \: : \: v|_f \in V^k(f) \textrm{ for all } f \textrm{ face of } \partial P
\Big\} .
$$
The above space is made of functions that on each face are two-dimensional virtual functions, that glue continuously across edges. Recalling Lemma \ref{lem:dofs1}, it follows that the following linear operators constitute a set of degrees of freedom for the space $\BdPk$:
\begin{align}
\label{Ldof1}
&\bullet\quad\textrm{evaluation of } v(\nu) \ \forall \nu \textrm{ vertex of } P ; \\
\label{Ldof2}
&\bullet\quad\textrm{evaluation of } v(\nu_e^i) \ \forall e \textrm{ edge of }\partial P, \ i=\{1,2,..,k-1 \} ; \\
\label{Ldof3}
&\bullet\quad\textrm{moments } \int_f v \: p_{k-2} \ \forall p_{k-2} \in \P_{k-2}(f),  \forall f \textrm{ face of } \partial P.
\end{align}
Once the boundary space is defined, the steps to follow in order to define the local virtual space on $P$ become very similar to the two dimensional case. We first introduce a preliminary local virtual element space on $P$
$$
\widetilde{V}^k(P) = \Big\{ v \in H^1(P) \: : \:  v|_{\partial P} \in \BdPk  , \
\Delta v \in \P_{k}(P) 
\Big\} 
$$
and the ``internal'' linear operators 
\begin{align}
\label{Ldof4}
&\bullet\quad\textrm{moments } \int_P v \: p_{k-2} \ \ \forall p_{k-2} \in \P_{k-2}(P)  .
\end{align}
We can now define the projection operator $\Pi^\nabla_P \: : \: \widetilde{V}^k(P) \rightarrow \P_k(P)$ by 
\begin{equation}
\left\{
\begin{aligned}
& \int_P \nabla (v - \Pi^\nabla_P  v)  \cdot \nabla p_k = 0 \quad \forall p_k \in \P_k(P) \\
& \text{for $k=1$:}\quad\sum_{\substack{\nu\,\text{vertex} \\[0.5mm]\text{ of}\,P}} \left(v(\nu) - \Pi^\nabla_P v(\nu)\right) = 0, \\
& \text{for $k\ge2$:}\quad\int_{ P} (v - \Pi^\nabla_P  v) = 0 .
\end{aligned}
\right.
\label{eqn:nablaProj}
\end{equation}

An integration by parts and observation \eqref{obs-1} show that the projection $\Pi^\nabla_P$ only depends on the operator values \eqref{Ldof1}-\eqref{Ldof2}-\eqref{Ldof3} and \eqref{Ldof4}. 
Therefore we can define the local virtual space
$$
V^k(P) = \Big\{ v \in \widetilde{V}^k(P) \: : \: \int_P v \, q = \int_P (\Pi^\nabla_P  v) q \textrm{ for all } q \in 
\P_{k-1}^{{\rm hom}}(P)\cup\P_{k}^{{\rm hom}}(P) 
\Big\} .
$$
The proof of the following lemma mimicks the two dimensional case, see for instance \cite{Ahmed-et-al:2013}.

\begin{lem}\label{lem:dofs2}
The operators \eqref{Ldof1}-\eqref{Ldof2}-\eqref{Ldof3} and \eqref{Ldof4} constitute a set of degrees of freedom for the space $V^k(P)$.
\end{lem}
It is immediate to verify that $\P_k(P) \subseteq V^k(P)$, that is a fundamental condition for the approximation properties of the space. Moreover, again due to the observation above, the projection operator 
$\Pi^\nabla_P : V^k(P) \rightarrow \P_k(P)$ is computable only on the basis of the degree of freedom values \eqref{Ldof1}-\eqref{Ldof2}-\eqref{Ldof3} and \eqref{Ldof4}. In addition, by following the same identical argument as in \eqref{obs-1} we obtain that $\int_P v p_{k}$ is computable for any $p_{k} \in \P_{k}(P)$ by using the degrees of freedom. Therefore also the $L^2$ projection operator $\Pi^0_P : V^k(P) \rightarrow \P_{k}(P)$, defined for any $v \in V^k(P)$ by
\begin{equation}
\int_P (v - \Pi^0_P  v)  q_k = 0 \quad \forall q_k \in \P_k(P) ,
\label{eqn:zeroProj}
\end{equation}
is computable by using the degree of freedom values.

Finally, the global virtual space $V^k \subset H^1(\Omega)$ is defined by using a standard assembly procedure as in finite elements. We define
$$
V^k = \Big\{ v \in H^1(\Omega) \: : \: v|_P \in V^k(P) \textrm{ for all } P \in \Omega_h \Big\} .
$$ 
The associated (global) degrees of freedom are the obvious counterpart of the local ones introduced above, i.e.
\begin{align}
\label{Gdof1}
&\bullet\quad\textrm{evaluation of } v(\nu) \ \forall \nu \textrm{ vertex of } \Omega_h\backslash\Gamma ; \\
\label{Gdof2}
&\bullet\quad\textrm{evaluation of } v(\nu_e^i) \ \forall e \textrm{ edge of } \Omega_h\backslash\Gamma, \ i=\{1,2,..,k-1 \} ; \\
\label{Gdof3}
&\bullet\quad\textrm{moments } \int_f v \: p_{k-2} \ \ \forall p_{k-2} \in \P_{k-2}(f), \forall f  \textrm{ face of } \Omega_h\backslash\Gamma; \\
\label{Gdof4}
&\bullet\quad\textrm{moments } \int_P v \: p_{k-2} \ \ \forall p_{k-2} \in \P_{k-2}(P), \forall P \in \Omega_h  .
\end{align}

\subsection{Discretization of the problem}\label{sub:probDisc}

We start by introducing the discrete counterpart of the involved bilinear forms. Given any polyhedron $P \in \Omega_h$ we need to approximate the local forms
$$
a_P(v,w)=\int_P \nabla v \cdot \nabla w \ , \quad m_P(v,w)=\int_P v \, w .
$$
We follow \cite{volley,Hitchhikers}. We first introduce the stabilization form
\begin{equation}\label{L:stab}
s_P(v,w) = \sum_{i=1}^{N^P_{\textrm{dof}}} \Xi_i(v) \: \Xi_i(w) \quad \forall v,w \in V^k(P) ,
\end{equation}
where $\Xi_i(v)$ is the operator that evaluates the function $v$ in the $i^{th}$ local degree of freedom and $N^P_{\textrm{dof}}$ denotes the number of such local degrees of freedom, see \eqref{Ldof1}-\eqref{Ldof2}-\eqref{Ldof3} and \eqref{Ldof4}.
We then set, for all $v,w \in V^k(P)$,
\begin{eqnarray}
a_P^h(v,w) &=& \int_P (\nabla \Pi^\nabla_P v) \cdot (\nabla \Pi^\nabla_P w) 
+ h_P \, s_P(v - \Pi^\nabla_P v , w - \Pi^\nabla_P w)  , \label{eqn:aForm}\\
m_P^h(v,w) &=& \int_P (\Pi^0_P v) (\Pi^0_P w) 
+ |P| \, s_P(v - \Pi^0_P v , w - \Pi^0_P w) \nonumber .
\end{eqnarray}
The above bilinear forms are consistent and stable in the sense of \cite{volley}. The global forms are given by, for all $v,w \in V^k$,
$$
a^h(v,w) = \sum_{P \in \Omega_h} a_P^h(v,w) \ , \quad
m^h(v,w) = \sum_{P \in \Omega_h} m_P^h(v,w) .
$$
Let now the discrete space with boundary conditions and its corresponding test space
$$
V^k_\Gamma = \Big\{ v \in V^k \: : \: v|_\Gamma = r_I \Big\} \ , \quad
V^k_0 = \Big\{ v \in V^k \: : \: v|_\Gamma = 0 \Big\} ,
$$
where $r_I$ is, face by face, an interpolation of $r$ in the virtual space $V^k(f)$.

We can finally state the discrete problem
\begin{equation}
\left\{
\begin{aligned}
& \textrm{Find } u_h \in V^k_\Gamma \textrm{ such that } \\
& a^h(u_h,v_h)  \: + \: m^h(u_h,v_h) \: = \: \int_\Omega f_h v_h + \int_{\Gamma'} g_h v_h
\qquad \forall v_h \in V^k_0 \ ,
\end{aligned}
\right.
\end{equation}
where the approximate loading $f_h$ is the $L^2$-projection of $f$ on piecewise polynomials of degree $k$, and where $g_h$ is the $L^2$-projection of $g$ on piecewise polynomials (still of degree $k$) living on $\Gamma'$. Note that all the forms and operators appearing above are computable in terms of the degree of freedom values of $u_h$ and $v_h$. 

We close this section by recalling a convergence result. The main argument for the proof can be found in \cite{volley}, while the associated interpolation estimates where shown in \cite{Steklov-VEM} for two dimensions and extended in \cite{Cangiani:apos} to the three dimensional case. 

Let now $\{ \Omega_h \}_h$ be a family of meshes, satisfying the following assumption. It exists a positive constant $\gamma$ such that all elements $P$ of $\{ \Omega_h \}_h$ and all faces of $\partial P$ are star-shaped with respect to a ball of radius bigger or equal than $\gamma h_P$; moreover all edges $e \in \partial P$, for all $P \in \{ \Omega_h \}_h$ have length bigger or equal than $\gamma h_P$.
\begin{thm}
Let the above mesh assumptions hold. Then, if the data and solution is sufficiently regular for the right hand side to make sense, it holds
\begin{equation}\label{H1-conv}
\| u - u_h \|_{H^1(\Omega)} \le C h^{s-1} \Big( |u|_{H^{s}(\Omega_h)} 
+ |f|_{H^{s-2}(\Omega_h)} + |g|_{H^{s-3/2}(\Gamma'_h)}
\Big) ,  
\end{equation}
where $2 \le s \le k+1$, the real $h$ denotes the maximum element diameter size and the constant $C$  is independent of the mesh size. The norms appearing in the right hand side are broken Sobolev norms with respect to the mesh (or its faces). 
\end{thm}

The above result applies also if $1 \le s<2$, but in that case the regularities on the data $f,g$ need to be modified. 
If the domain $\Omega$ is convex (or regular) then under the same assumptions and notations it also holds \cite{variable-primale}
\begin{equation}\label{L2-conv}
\| u - u_h \|_{L^2(\Omega)} \le C h^{s} \Big( |u|_{H^{s}(\Omega_h)} 
+ |f|_{H^{s-2}(\Omega_h)} + |g|_{H^{s-3/2}(\Gamma'_h)}
\Big) . 
\end{equation}
Finally, we note that an extension of the results to more general mesh assumptions could be possibly derived following the arguments for the two-dimensional case shown in \cite{stab:theory}. 

In the following numerical tests, we will in particular show the robustness of the method also for quite irregular meshes.

\section{Numerical tests}\label{NUMS}

In this section we collect the numerical results to evaluate 
the reliability and robustness of the Virtual Element Method in three dimensions.

\subsection{Meshes and error estimators}\label{ssec:mesh}

Before dealing with the numerical examples, 
we define the domains where we solve the PDEs and 
the polyhedral meshes which discretize such domains.
Moreover, we define the norms that we use
to evaluate the error.

\subsubsection{Meshes}\label{sub::mesh}

\begin{figure}[!htb]
\begin{center}
\begin{tabular}{ccc}
\includegraphics[width=0.26\textwidth]{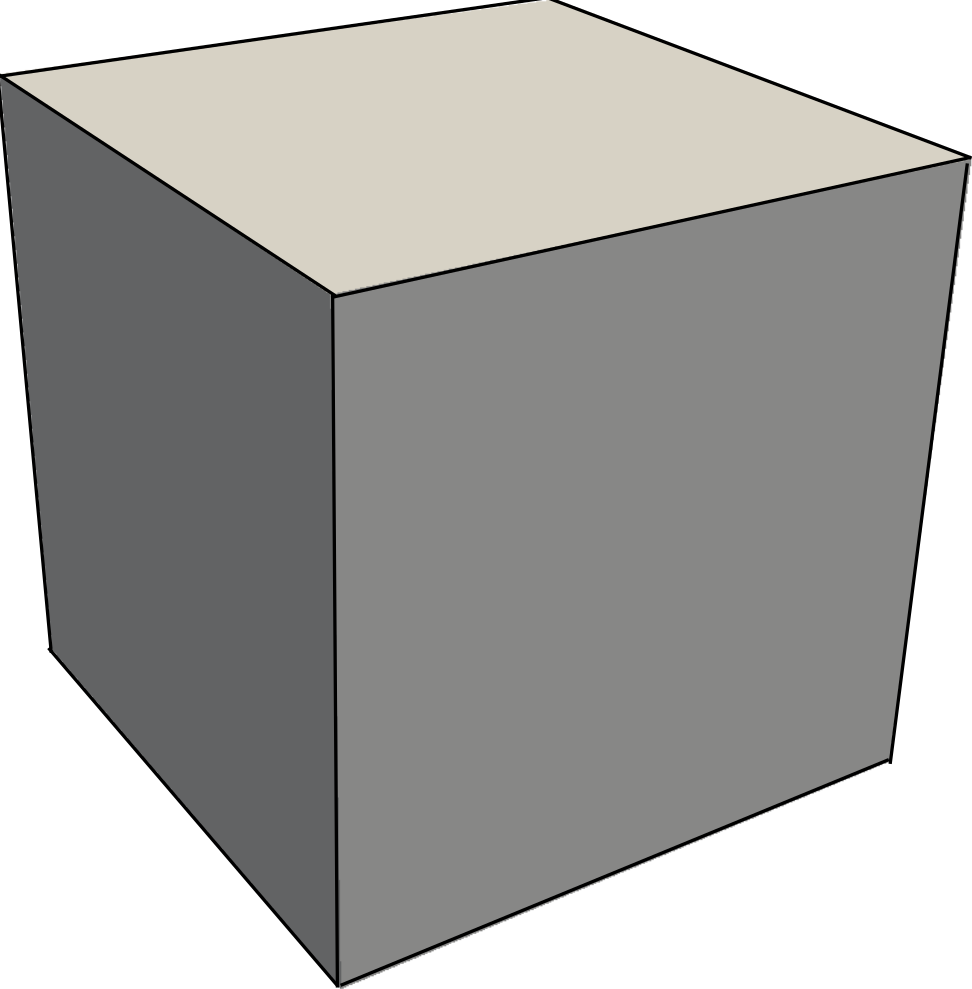} &\phantom{mm} &
\includegraphics[width=0.26\textwidth]{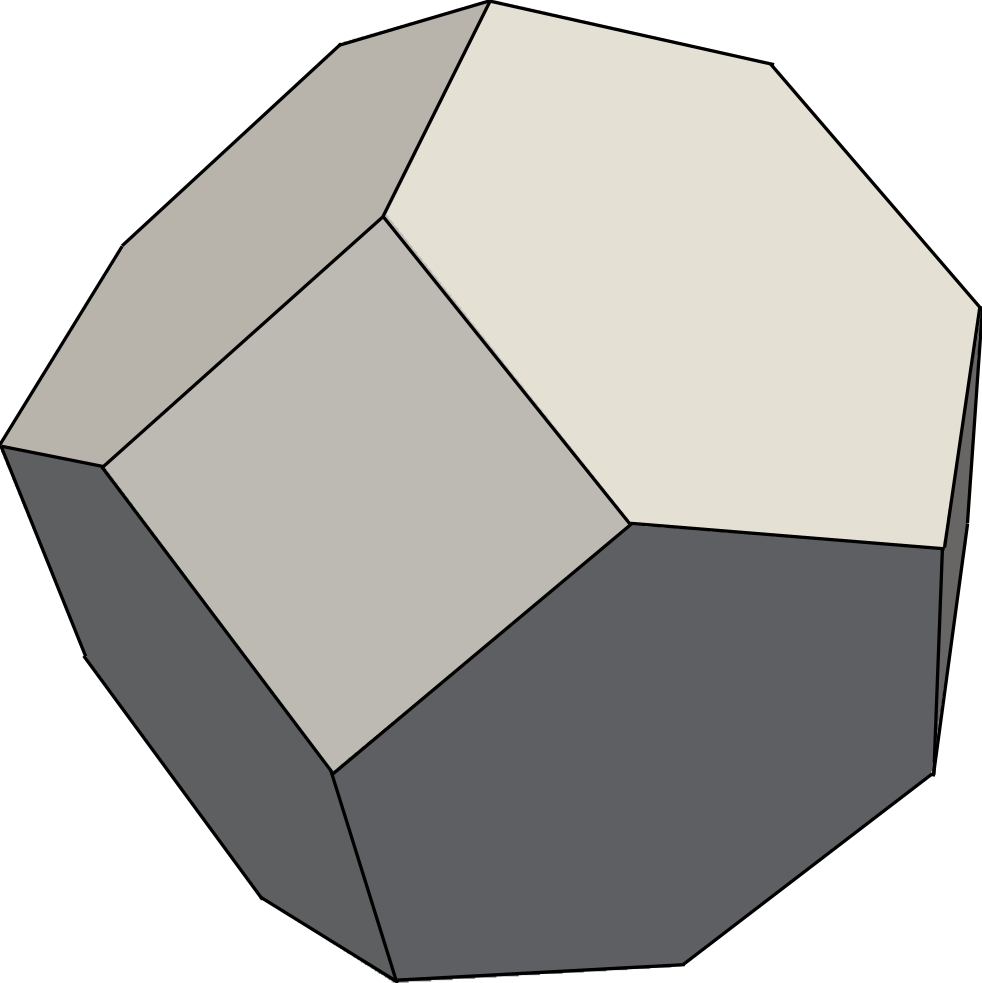} \\
(a) & & (b)
\end{tabular}
\end{center}
\caption{(a) Standard cube $[0,\,1]^3$ and (b) the truncated octahedron.}
\label{fig:geo}
\end{figure}

We consider two different domains: 
the standard $[0,\,1]^3$ cube, see Figure~\ref{fig:geo} (a),
and a truncated octahedron~\cite{patagonGeo}, see Figure~\ref{fig:geo} (b).

We make different discretizations of such domains by exploiting the c++ library \texttt{voro++}~\cite{voroPlusPlus}.
More specifically we will consider the following three mesh types:

\phantom{a}

\noindent \textbf{Random} refers to a mesh
where the control points of the Voronoi tessellation are randomly displaced inside the domain.
We underline that these kind of meshes are characterized by stretched polyhedrons so 
the robustness of the VEM will be severely tested.

\phantom{a}

\noindent \textbf{CVT} refers to a Centroidal Voronoi Tessellation, i.e., 
a Voronoi tessellation where the control points coincides with the centroid of the cells they define.
We generate such meshes via a standard Lloyd algorithm~\cite{cvtPaper}.
In this case the Voronoi cells are more regular than the ones of the previous case.

\phantom{a}

\begin{figure}[!htb]
\begin{center}
\begin{tabular}{ccc}
\includegraphics[width=0.35\textwidth]{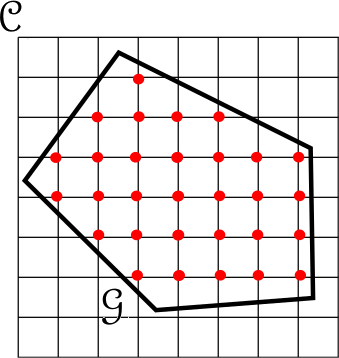} &\phantom{mm} & 
\includegraphics[width=0.35\textwidth]{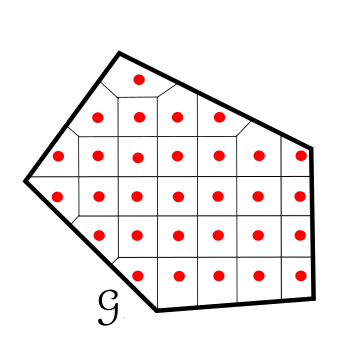} \\
(a) & & (b)
\end{tabular}
\end{center}
\caption{(a) The structured mesh of a square $\mathcal{C}$ which contains the input geometry $\mathcal{G}$,
where we highlight the vertices inside $\mathcal{G}$. 
(b) The final polyhedral Voronoi mesh where we highlight the control points.}
\label{fig:constr}                                                                       
\end{figure}

\noindent \textbf{Structured} refers to meshes composed by structured cubes inside the domain and 
arbitrary shaped mesh close to the boundary.
This mesh is built by considering as control points the vertices of a structured mesh of a cube $\mathcal{C}$ 
(containing the input geometry $\mathcal{G}$) which are inside $\mathcal{G}$, 
see Figure~\ref{fig:constr} for a two dimensional example.
When we consider the cube $[0,1]^3$, 
this kind of mesh coincides with a structured mesh composed by cubes.
These meshes are really interesting from the computational point of view.
Indeed, the VEM local matrices are exactly the same for all cubes inside the domain.
It is therefore possible to compute such local matrices \emph{only} once so
the computational effort in assembling the stiffness matrix as well as the right hand side is reduced.
Moreover, the ensuring scheme may inherit some advantages of structural cubic meshes.

\begin{figure}[!htb]
\begin{center}
\begin{tabular}{ccc}
\multicolumn{3}{c}{{\large Random}}\\
\includegraphics[width=0.4\textwidth]{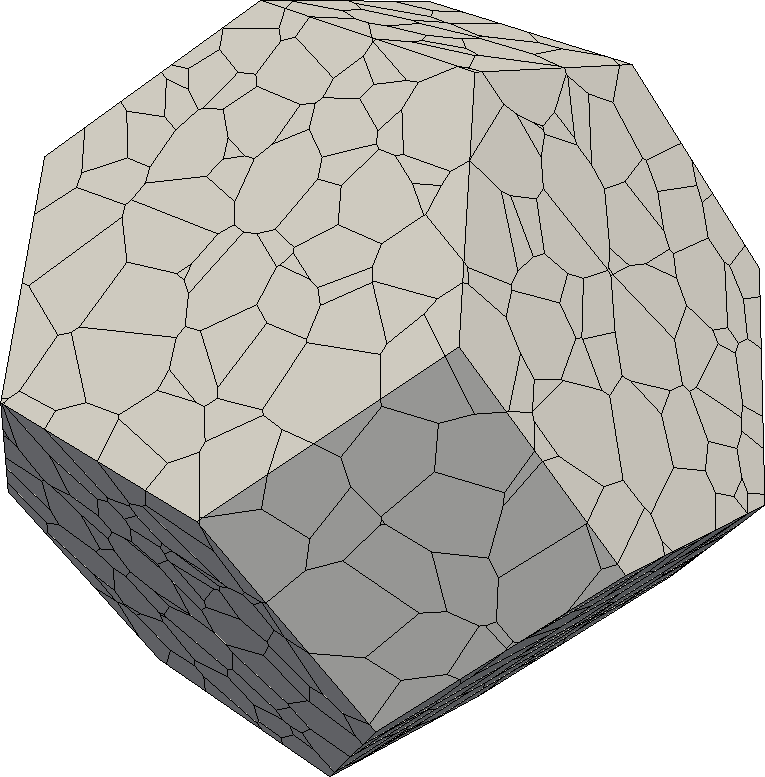} & &\phantom{mm}
\includegraphics[width=0.4\textwidth]{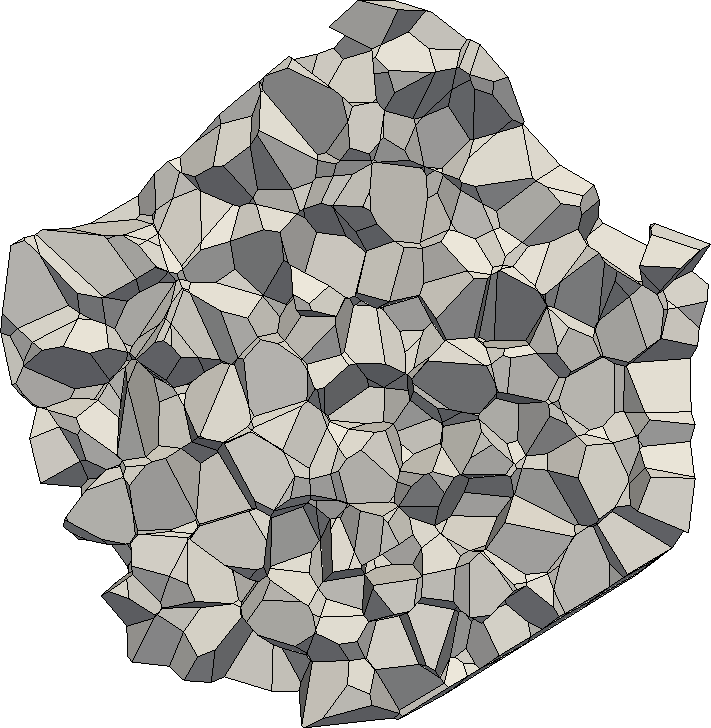} \\
\multicolumn{3}{c}{{\large CVT}}\\
\includegraphics[width=0.4\textwidth]{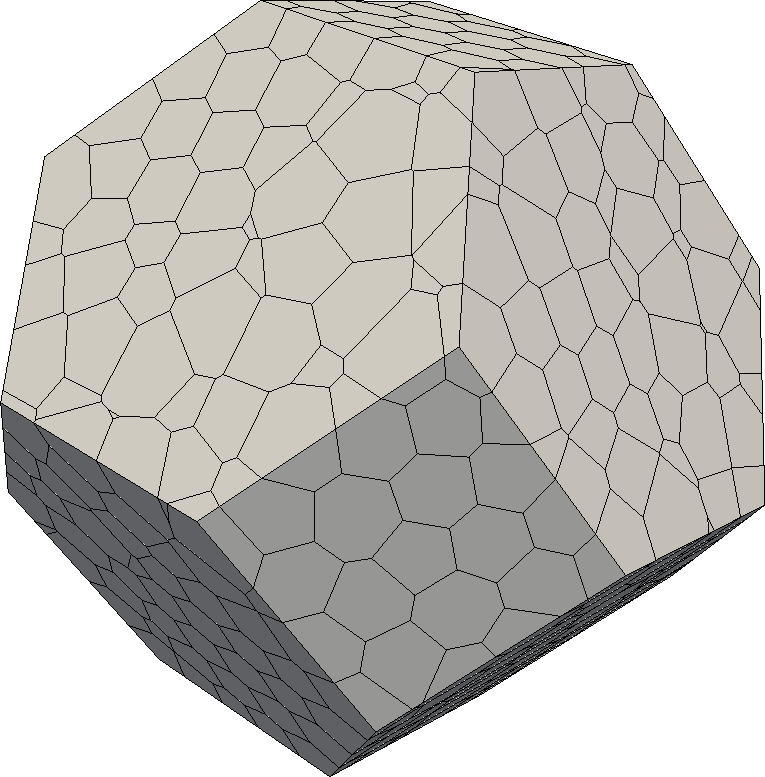} & &\phantom{mm}
\includegraphics[width=0.4\textwidth]{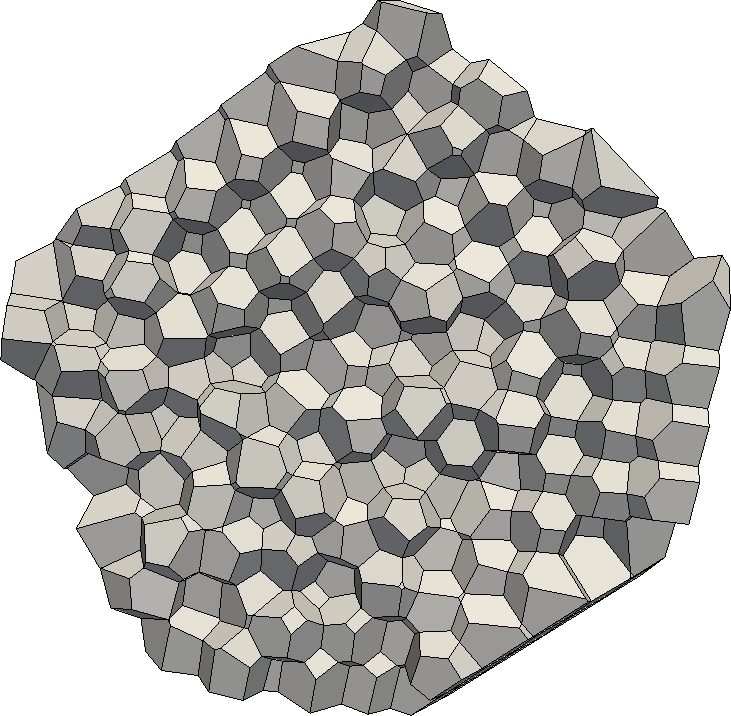} \\
\multicolumn{3}{c}{{\large Structured}}\\
\includegraphics[width=0.4\textwidth]{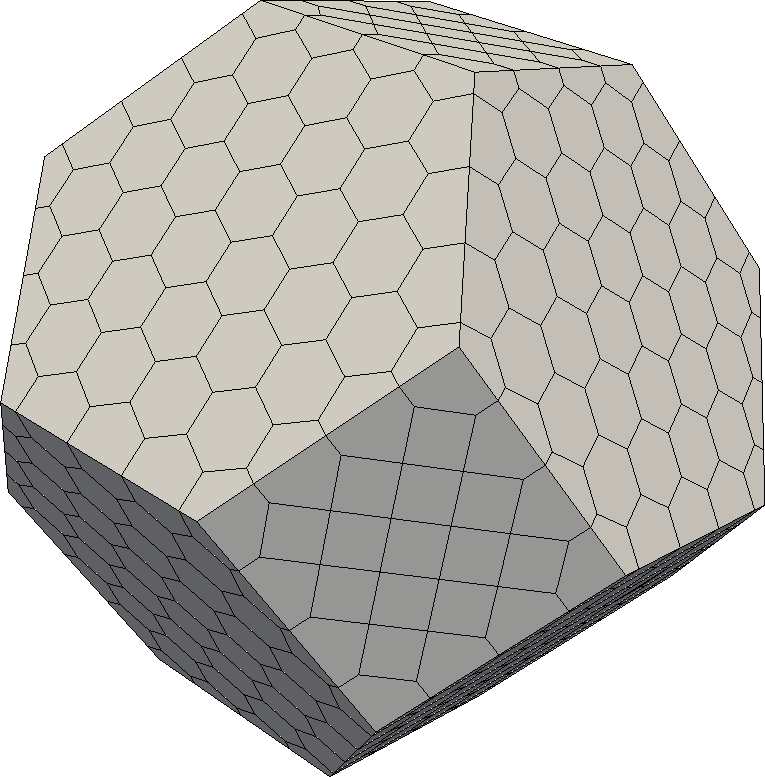} & &\phantom{mm}
\includegraphics[width=0.4\textwidth]{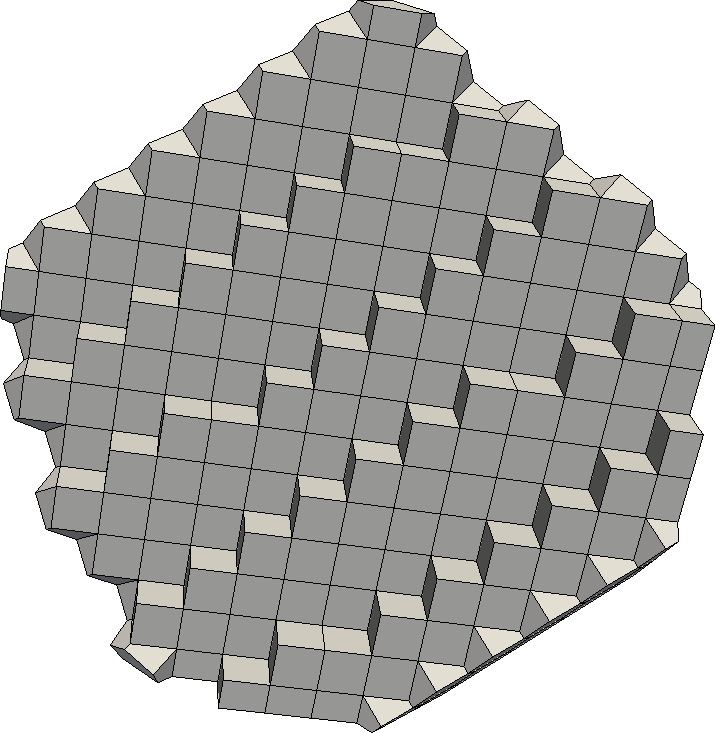} \\
\end{tabular}
\end{center}
\caption{Different discretization of the truncated octahedron.}
\label{fig:meshes}
\end{figure}

\medskip
In Figure~\ref{fig:meshes} we collect an example on the truncated octahedron geometry of all these kinds of meshes.
To analyze the error convergence rate,
we make a sequence of meshes with decreasing size for each mesh type.

\subsubsection{Error norms}

Let $u$ be the exact solution of the PDE and $u_h$ the discrete solution provided by the VEM.
To evaluate how this discrete solution is close to the exact one,
we use the local projectors of degree $k$ on each polyhedron $P$ of the mesh, 
$\Pi^{\nabla}_P\,u_h$ and $\Pi^0_P\,u_h$, 
defined in~\eqref{eqn:nablaProj} and~\eqref{eqn:zeroProj}, respectively.
We compute the following quantities:
\begin{itemize}
\item\textbf{$\mathbf{H^1}$-seminorm error}
$$
e_{H^1} := \sqrt{\sum_{P\in\Omega_h}\left| u - \Pi^{\nabla}_P\,u_h \right|^2_{H^1(P)}}\,,
$$
\item\textbf{$\mathbf{L^2}$-norm error}
$$
e_{L^2} := \sqrt{\sum_{P\in\Omega_h}\left\| u - \Pi^0_P\,u_h \right\|^2_{L^2(P)}}\,,
$$
\item\textbf{$\mathbf{L^{\infty}}$-norm error}
$$
e_{L^\infty} := \max_{\nu\in\mathcal{N}}\left|\,u(\nu) - u_h(\nu)\,\right|\,,
$$
where $\mathcal{N}$ is the set of all the vertexes and internal edge nodes of the VEM scheme,
see Equations~\eqref{Gdof1} and~\eqref{Gdof2}.
Since we do not take the $\max$ over all the domain but \emph{only} on some nodes,
$e_{L^\infty}$ is an approximation of the true $L^\infty$-norm.
Moreover, in this case we can directly compute such quantity without resorting to the projections operators.
\end{itemize}

In the following subsections we will present some numerical tests 
to underline different computational aspects of the method.
In all cases the mesh-size parameter $h$ is measured in an averaged sense 
\begin{equation}\label{h-medio}
h = \left( \frac{|\Omega|}{N_P} \right)^{1/3} ,
\end{equation}
with $N_P$ denoting the number of polyhedrons in the mesh.

\subsection{Test case 1: $h$-analysis for diffusion problem on a cube}

Let us consider the problem 
\begin{equation}
\left\{
\begin{array}{rl}
-\Delta u &=\, f\quad\quad\textnormal{in }\Omega\\
u &=\,r \quad\quad\textnormal{on }\Gamma\\
\frac{\partial u}{\partial n} &=\,0  \quad\quad\textnormal{on }\partial\Omega\backslash\Gamma 
\end{array}
\right.,
\label{eqn:neudiri}
\end{equation}
where the domain $\Omega$ is the cube $[0,\,1]^3$ and 
$\Gamma$ is the union of the four faces corresponding to the planes $y=0$, $y=1$, $z=0$ and $z=1$.
We choose the right hand side $f$ and $r$ in such a way that the exact solution is 
$$
u(x,\,y,\,z):=\sin(\pi x)\cos(\pi y)\cos(\pi z)\,.
$$
In this example we consider all the three types of discretizations introduced in Subsection~\ref{ssec:mesh}, 
i.e., Random, CVT and Structured. 
Note that in this case the last type of discretization becomes a standard structured cubic mesh.

\begin{figure}[!htb]
\begin{center}
\begin{tabular}{cc}
\includegraphics[width=0.49\textwidth]{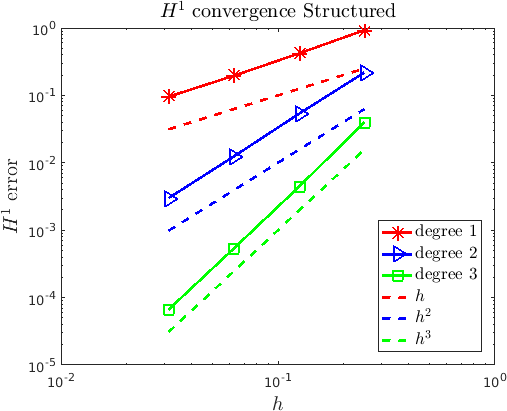} &
\includegraphics[width=0.49\textwidth]{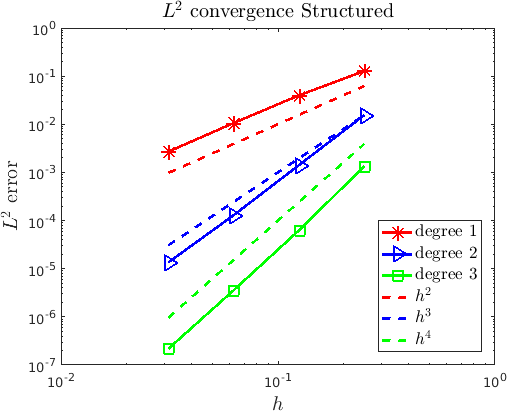} \\
\includegraphics[width=0.49\textwidth]{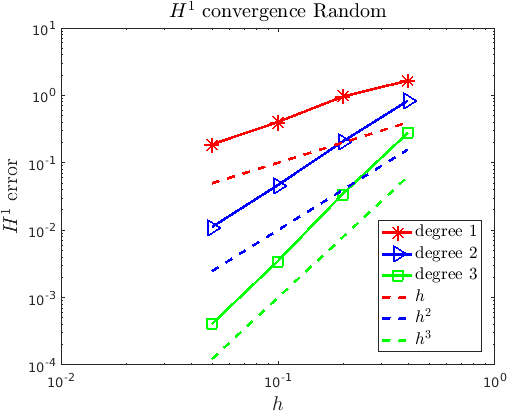} &
\includegraphics[width=0.49\textwidth]{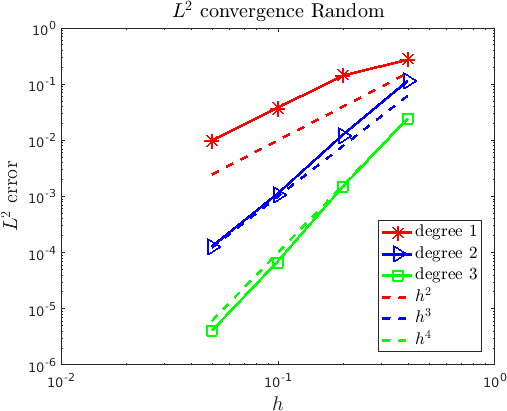} \\
\includegraphics[width=0.49\textwidth]{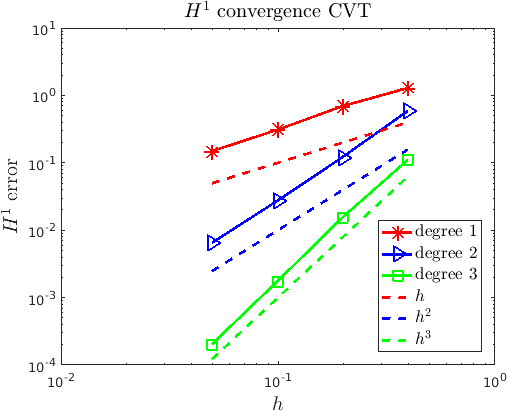} &
\includegraphics[width=0.49\textwidth]{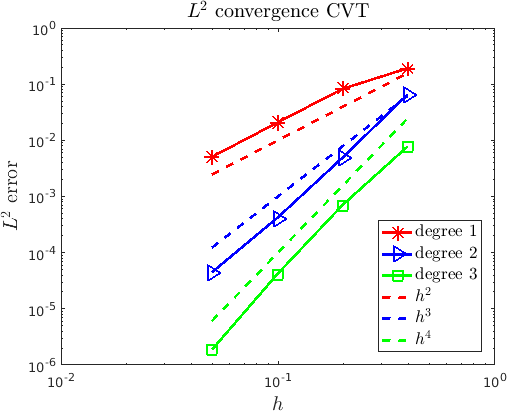} \\
\end{tabular}
\caption{Test case 1: $h$-convergence with different meshes.}
\label{fig:convNeuDiri}
\end{center}
\end{figure}

In Figure~\ref{fig:convNeuDiri} we show the resulting graphs and 
in Tables~\ref{tab:convNeuDiriH1} and~\ref{tab:convNeuDiriL2} we provide the convergence rates.
These data show that we achieve the theoretical convergence rate for all the VEM approximation degrees and 
for each type of meshes, see equations \eqref{H1-conv}, \eqref{L2-conv}.

\begin{table}[!htb]
\begin{center}
\begin{tabular}{|l|l|l|l|}
\cline{2-4}
\multicolumn{1}{}{}&\multicolumn{3}{|c|}{convergence rates} \\
\hline
mesh type &$k=1$ &$k=2$&$k=3$ \\
\hline
Structured  &1.0344 &2.0543 &3.0125 \\
Random      &1.0927 &2.0465 &3.0656 \\
CVT         &1.0492 &2.0672 &3.0970 \\
\hline
\end{tabular}
\end{center}
\caption{Test case 1: $H^1$-seminorm convergence rates.}
\label{tab:convNeuDiriH1}
\end{table}

\begin{table}[!htb]
\begin{center}
\begin{tabular}{|l|l|l|l|}
\cline{2-4}
\multicolumn{1}{}{}&\multicolumn{3}{|c|}{convergence rates} \\
\hline
mesh type &$k=1$ &$k=2$&$k=3$ \\
\hline
Structured  &1.9763 &3.2551 &4.0372 \\
Random      &1.9136 &3.1067 &4.0678 \\
CVT         &2.0230 &3.2144 &4.4620 \\
\hline
\end{tabular}
\end{center}
\caption{Test case 1: $L^2$-norm convergence rates.}
\label{tab:convNeuDiriL2}
\end{table}

\subsection{Test case 2: $h$-analysis for diffusion-reaction problem on a polyhedron}
In this example we consider the problem 
\begin{equation}
\left\{
\begin{array}{rl}
-\Delta u + u &=\, f\quad\quad\textnormal{in }\Omega\\
u &=\,r \quad\quad\textnormal{on }\Gamma = \partial\Omega\\
\end{array}
\right.,
\label{eqn:reacDiff}
\end{equation}
where the domain $\Omega$ is the truncated octahedron, see Figure~\ref{fig:geo} (b),
and we choose the right hand side $f$ and $r$ in such a way that the exact solution is 
$$
u(x,\,y,\,z):=\sin(2xy)\,\cos(z)\,.
$$
In this example we analyze the convergence rate for VEM approximation degrees from 1 to 3,
and compare the results obtained with the three mesh types described in Subsection~\ref{sub::mesh}.

\begin{figure}[!htb]
\begin{center}
\begin{tabular}{cc}
\includegraphics[width=0.49\textwidth]{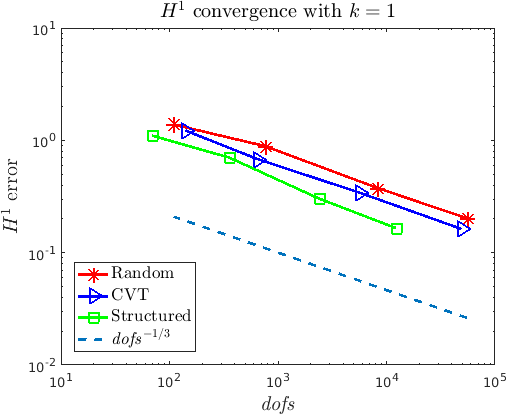} &
\includegraphics[width=0.49\textwidth]{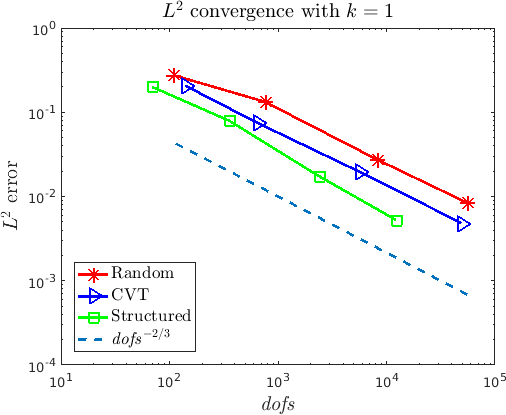} \\
\includegraphics[width=0.49\textwidth]{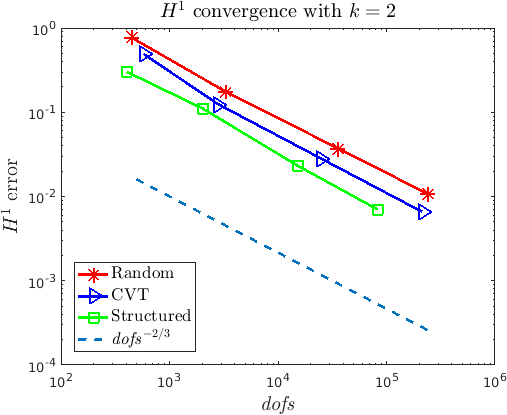} &
\includegraphics[width=0.49\textwidth]{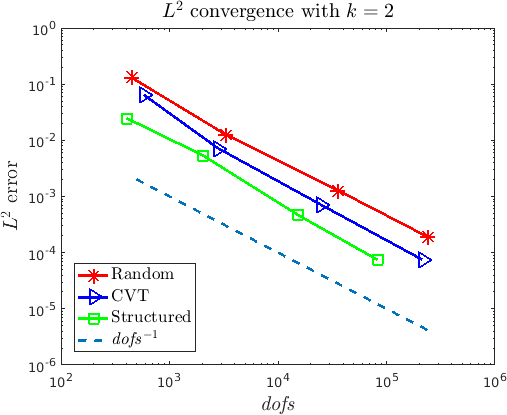} \\
\includegraphics[width=0.49\textwidth]{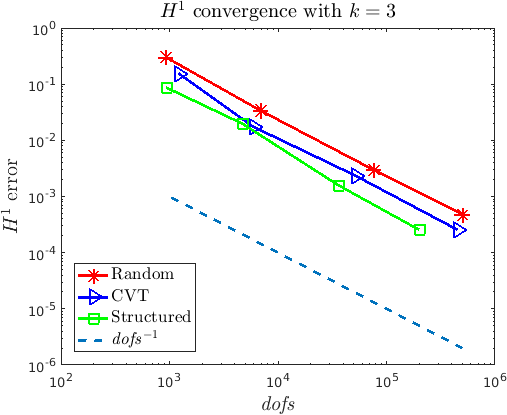} &
\includegraphics[width=0.49\textwidth]{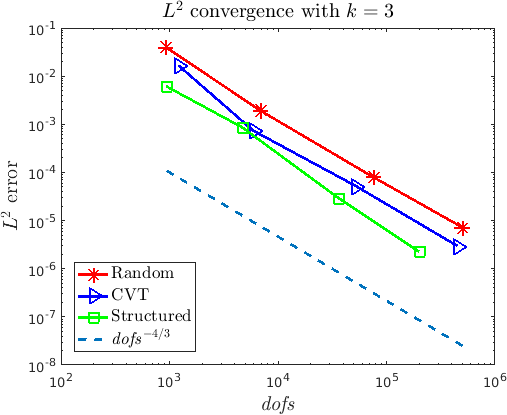} \\
\end{tabular}
\caption{Test case 2: $dofs$-convergence with different meshes.}
\label{fig:convReacDiff}
\end{center}
\end{figure}

In Figure~\ref{fig:convReacDiff} we plot the convergence graphs 
with respect to the total number of degrees of freedom $N_{\textrm{dof}}$.  
Considering that (for fixed order $k$ and for fixed mesh family) the mesh-size parameter is expected to behave as
$
h \sim N_{\textrm{dof}}^{-1/3},
$
it follows that both error norms behave as expected from the theory (see \eqref{H1-conv} and \eqref{L2-conv}).

Moreover we observe that the error is slightly affected by the shape of the mesh elements.
Indeed, the errors associated with the Random mesh are always larger than the ones obtained with
a more regular mesh, while the Structured meshes yield the best results (even when compared to the CVT meshes).

\subsection{Test case 3: convergence analysis with different $k$}

In this example we consider the convergence with respect to the accuracy degree $k$.
We fix the truncated octahedron geometry and we solve the following problem
\begin{equation}
\left\{
\begin{array}{rl}
-\Delta u &=\, f\quad\quad\textnormal{in }\Omega\\
u &=\,r \quad\quad\textnormal{on }\Gamma=\partial \Omega\\
\end{array}
\right.,
\label{eqn:pconvergence}
\end{equation}
where the right hand side $f$ and the boundary condition $r$ are chosen in such a way that
the exact solution is
$$
u(x,\,y,\,z):= \sin(\pi x)\,\cos(\pi y)\,\cos(\pi z)\,.
$$

\begin{figure}[!htb]
\begin{center}
\begin{tabular}{cc}
\includegraphics[width=0.47\textwidth]{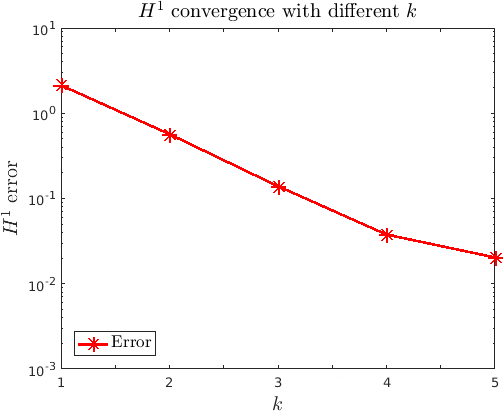} &
\includegraphics[width=0.47\textwidth]{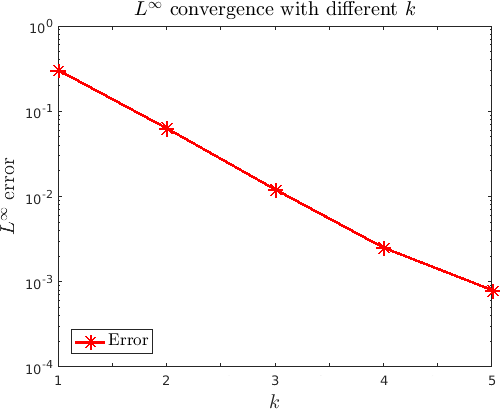} \\
\end{tabular}
\caption{Test case 3: convergence with respect to the order $k$.}
\label{fig:pconvergenceNoStabCase5}
\end{center}
\end{figure}

The mesh is kept fixed (the CVT mesh of the truncated octahedron composed by 116 polyhedrons)
and we rise the polynomial degree $k$ from 1 to 5.
In Figure~\ref{fig:pconvergenceNoStabCase5} we provide the convergence graphs of both $H^1$-seminorm and $L^\infty$-norm.
The trend of these errors show an exponential convergence in terms of $k$ and are thus
aligned with the existing theory for the two-dimensional case \cite{hp-uniform}.
However, both the $H^1$ and the $L^\infty$ errors show a slight bend in the convergence graphs for $k=5$.
This behavior is probably due to the stabilizing matrix~\eqref{L:stab}
that should be better devised in order to develop a spectral approximation strategy. 

Although such aspect deserves a deeper study, that is beyond the scopes of the current paper, 
we propose a novel stabilization strategy which, at least in the present context, cures the problem. 

Let $\{ \varphi_i \}_{i=1}^{N^P_{{\rm dof}}}$ represents the canonical basis functions on element $P$,
defined by
$$
\varphi_i \in V^k(P) \ , \qquad \Xi_j (\varphi_i) = \delta_{ij} \ \textrm{ for } \ j=1,2,..,N^P_{{\rm dof}},
$$
where we refer to \eqref{L:stab} for the notation.
Then, a diagonal stabilization form $s_P(\cdot,\cdot)$ should satisfy (see \eqref{eqn:aForm} and \cite{volley,Hitchhikers})
$$
h_P s_P(\varphi_i,\varphi_i) \simeq a_P(\varphi_i,\varphi_i) \qquad i=1,2,..,N^P_{{\rm dof}}
$$
in order to mimic the original energy $a_P(\cdot,\cdot)$ of the basis functions. The original form in \eqref{L:stab} corresponds to assuming $a_P(\varphi_i,\varphi_i) \simeq h_P$, that (considering the involved scalings) can be shown to be a reasonable choice with respect to $h_P$ (and thus for moderate $k$). On the other hand, choice \eqref{L:stab} may be less suitable for a higher $k$, especially in three dimensions, where different basis functions may carry very different energies. We therefore propose a simple alternative (still using a diagonal stabilizing form)
$$
\widetilde{s}_P(\varphi_i,\varphi_i) = \max \{ h_P , a_P(\Pi^\nabla_P\varphi_i,\Pi^\nabla_P\varphi_i) \} \qquad i=1,2,..,N^P_{{\rm dof}} .
$$
Note that at the practical level computing $a_P(\Pi^\nabla_P\varphi_i,\Pi^\nabla_P\varphi_i)$ is immediate, since such term is simply the $i^{th}$ term on the diagonal of the consistency matrix $ a_P ( \Pi^\nabla_P \varphi_i , \Pi^\nabla_P \varphi_j) $, that is already computed.

\begin{figure}[!htb]
\begin{center}
\begin{tabular}{cc}
\includegraphics[width=0.47\textwidth]{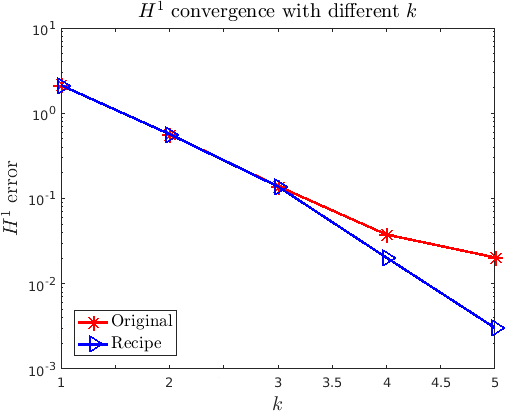} &
\includegraphics[width=0.47\textwidth]{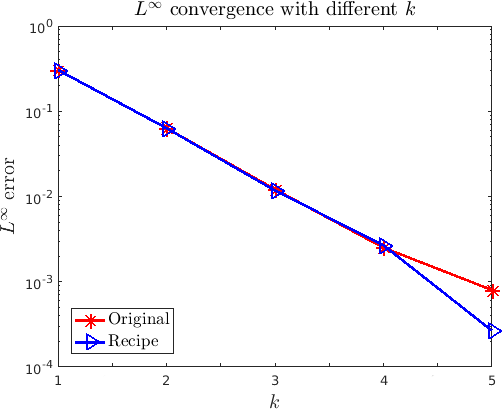} \\
\end{tabular}
\caption{Test case 3: convergence with respect to the order $k$, 
comparison between the original and the diagonal recipe stabilizations.}
\label{fig:recipe}
\end{center}
\end{figure}

The choice above, referred in the following as \emph{diagonal recipe}, cures the problem, see Figure~\ref{fig:recipe}. 
This fact is further confirmed by the data in Table~\ref{tab:recipePendence}.
Here, we compute the slopes of these lines at each step and 
we numerically show that the recipe diagonal stabilization yields better results.

\begin{table}[!htb]
\begin{center}
\begin{tabular}{|l|l|l|l|l|}
\cline{2-5}
\multicolumn{1}{c}{}&\multicolumn{4}{|c|}{$H^1$ convergence}\\
\hline{2-5}
Original &-1.3219 &-1.4065 &-1.3045 &-0.6195\\
\hline
Recipe   &-1.3218 &-1.4144 &-1.9211 &-1.8977\\
\hline
\multicolumn{5}{c}{}\\
\multicolumn{5}{c}{}\\
\cline{2-5}
\multicolumn{1}{c}{}&\multicolumn{4}{|c|}{$L^\infty$ convergence}\\
\hline
Original &-1.5704  &-1.6532 &-1.5736 &-1.1455\\
\hline
Recipe   &-1.5698  &-1.6850 &-1.4805 &-2.3065\\
\hline
\end{tabular}
\end{center}
\caption{Test case 3: slopes at each step for the original and the recipe diagonal stabilization.}
\label{tab:recipePendence}
\end{table}

\subsection{Test case 4: patch test}

In~\cite{volley,Hitchhikers} it is shown that VEM passes the so-called ``patch test''.
If we are dealing with a PDE whose solution is a polynomial of degree $k$ and 
we use a VEM approximation degree equal to $k$, we recover the ``exact solution'', i.e., 
the solution up to the machine precision.

We make a patch test for VEM approximation degrees $k$ from 1 up to 5.
More specifically, we consider the following PDE 
\begin{equation}
\left\{
\begin{array}{rl}
-\Delta u &=\, f\quad\quad\textnormal{in }\Omega\\
u &=\,r \quad\quad\textnormal{on }\Gamma=\partial \Omega\\
\end{array}
\right.
\label{eqn:patchTest}
\end{equation}
where the right hand side $f$ and the Dirichlet boundary condition $r$ are chosen
in accordance with the exact solution
$$
u(x,\,y,\,z):= (x+y+z)^k\,.
$$
Since we are not interested in varying the mesh size but \emph{only} the VEM approximation degree,
we run the experiments on the same coarse CVT mesh of the truncated octahedron composed by 116 polyhedrons.

In Table~\ref{tab:patch} we collect the results.
The errors are close to the machine precision,
but for higher VEM approximation degrees they become larger.
This fact is natural and stems from the conditioning of the matrices involved in the computation of the VEM solution.
Indeed, as in standard FEM, 
their condition numbers become larger when we consider higher VEM approximation degrees.

\begin{table}[!htb]
\begin{center}
\begin{tabular}{|c|c|c|c|}
\hline
solution degree &$H^1$-seminorm &$L^2$-seminorm &$L^\infty$-norm\\
\hline
1 &5.9775e-12 &6.1919e-13 &1.1479e-12\\                      
2 &2.2008e-11 &1.8416e-12 &4.9409e-12\\
3 &1.0490e-10 &1.1284e-11 &8.5904e-12\\
4 &3.0959e-10 &1.0039e-10 &2.6197e-11\\
5 &1.1563e-09 &6.9693e-09 &1.5433e-10\\
\hline
\end{tabular}
\end{center}
\caption{Test case 4: patch tests errors with the mesh CVT of the truncated octahedron composed by 116 elements.}
\label{tab:patch}
\end{table}

\subsection{Test case 5: stabilizing parameter analysis}

In this example we make an analysis on the stabilizing part of the local stiffness matrix.
We slightly modify $a_P^h$ defined in Subsection~\ref{sub:probDisc} 
by introducing the parameter $\tau\in\mathbb{R}$, i.e.,
\begin{equation}
a_P^h(v,w) = \int_P (\nabla \Pi^\nabla_P v) \cdot (\nabla \Pi^\nabla_P w)
+ \tau\,h_P \, s_P(v - \Pi^\nabla_P v , w - \Pi^\nabla_P w)  . 
\label{eqn:aFormWithTau}
\end{equation}
We fix a standard Poisson problem 
\begin{equation}
\left\{
\begin{array}{rl}
-\Delta u &=\, f\quad\quad\textnormal{in }\Omega\\
u &=\,0 \quad\quad\textnormal{on }\Gamma=\partial \Omega\\
\end{array}
\right.,
\label{eqn:stabTest}
\end{equation}
where $\Omega=[0,\,1]^3$ and whose exact solution is 
$$
u(x,\,y,\,z):= \sin(\pi x)\,\cos(\pi y)\,\cos(\pi z)\,.
$$
We exploit the same discretization of $\Omega$, a CVT mesh composed by 1024 polyhedrons,
and solve this problem for many values of the parameter $\tau$.
More specifically, we consider $\tau = 10^{t}$ for one hundred uniformly distributed values of $t$ in $[-2,2]$ and
we compute the errors for the corresponding VEM solutions in the $H^1$-seminorm and $L^\infty$-norm.

\begin{figure}[!htb]
\begin{center}
\begin{tabular}{cc}
\includegraphics[width=0.47\textwidth]{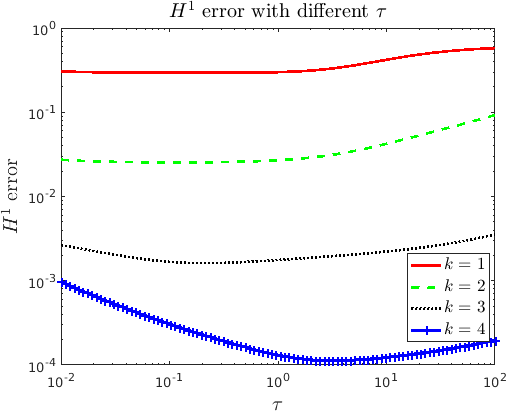} &
\includegraphics[width=0.47\textwidth]{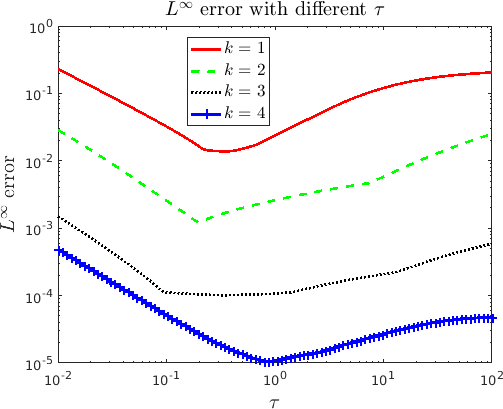} \\
(a) & (b)
\end{tabular}
\caption{Test case 5: (a) $H^1$ error and (b) $L^\infty$ error for varying $\tau$ and different 
VEM approximation degrees $k$.}
\label{fig:tauAnalysis}
\end{center}
\end{figure}

In Figure~\ref{fig:tauAnalysis} we provide the graphs of both errors, as a function of $\tau$ in a logarithmic scale, for VEM approximation degree $k=1,2,3$ and $4$.
We can observe that the trend of the errors is the same for each $k$: it grows whenever very small or very large choices are taken for the $\tau$ parameter. 
Although the $L^\infty$-norm seems more sensible than the $H^1$-seminorm, the method appears in general quite robust with respect to the parameter choice. For instance, we show in Table \ref{tab:ratios} the quantities
$$
\delta_{e_{H^1}} := \frac{\max_{\tau \in [10^{-1},10]}\left(e_{H_1}\right)}{\min_{\tau\in [10^{-1},10]}\left(e_{H_1}\right)} \ , \qquad
\delta_{e_{L^\infty}} := \frac{\max_{\tau \in [10^{-1},10]}\left(e_{L^\infty}\right)}{\min_{\tau\in [10^{-1},10]}\left(e_{L^\infty}\right)} \ , 
$$
representing the ratios of maximum to minimum error for the parameter range $ [10^{-1},10]$, that corresponds to a factor of 100 between minimum and maximum $\tau$.
From this table we can appreciate that the error ratios are all within an acceptable range.

\begin{table}[!htb]
\begin{center}
\begin{tabular}{|c|c|c|c|}
\hline
$k$  &$\delta_{e_{H^1}}$ &$\delta_{e_{L^\infty}}$ \\
\hline
1 &1.4487e+00 &8.9624e+00 \\                      
2 &1.7014e+00 &4.9538e+00 \\
3 &1.3884e+00 &2.0908e+00 \\
4 &2.6989e+00 &4.8393e+00 \\
\hline
\end{tabular}
\end{center}
\caption{Test case 5: Ratios of the errors for different degrees $k$. Parameter $\tau \in [10^{-1},10]$.}
\label{tab:ratios}
\end{table}


\begin{center} {\bf Aknowledgments} \end{center}
The first and second authors have received funding from the European Research Council (ERC) under the European Unions Horizon 2020 research and innovation programme (grant agreement no. 681162).

\bibliography{VEM}


\end{document}